# Maximizing General Set Functions by Submodular Decomposition

Kevin Byrnes
The Johns Hopkins University
Department of Applied Mathematics and Statistics
kbyrnes2@jhu.edu

**Abstract:** We present a branch and bound method for maximizing an arbitrary set function $\theta : 2^V \to \mathbb{R}$. By decomposing $\theta$ as f-$\delta$, where f is a submodular function and $\delta$ is the cut function of a (simple, undirected) graph G with vertex set V, our original problem is reduced to a sequence of submodular maximization problems. We characterize a class of submodular functions, which when maximized in the subproblems, lead the algorithm to converge to a global maximizer of f-$\delta$. Two "natural" members of this class are analyzed; the first yields polynomially-solvable subproblems, the second, which requires less branching, yields NP-hard subproblems but is amenable to a polynomial-time approximation algorithm. These results are extended to problems where the solution is constrained to be a member of a subset system. Structural properties of the maximizer of f-$\delta$ are also proved.

## §1. Introduction

Submodular functions arise naturally in many areas of combinatorial optimization. For example: the maximization of influence in a social network [11], optimal placement of sensors [12], and feature selection in machine learning [10]. In this paper we propose an additional use of submodular functions, to provide a "universal" decomposition of set functions $\theta : 2^V \to \mathbb{R}$. This decomposition will permit the maximization of *any* set function by a novel branch and bound procedure, where each subproblem consists of maximizing the decomposed function over a family of independent sets of a specially structured graph. After relaxation, each subproblem becomes equivalent to maximizing a certain submodular function.

While maximizing a submodular function is, in general, NP-hard, they possess significant exploitable structure. In this case, we investigate two methods of relaxation to the maximization of a submodular function: the first permits polynomially-solvable subproblems, while the second (tighter) method yields a polynomial-time constant-factor approximation algorithm [4] if our original function $\theta$ is non-negative. Both relaxations "work" in the context of our branch and bound framework in that their (exact or approximate) solutions permit valid pruning and fathoming rules. Variations of the original problem in which the solution is constrained to lie within a *subset system* are also explored, and our branch and bound and relaxation-solution techniques are extended to this case. Finally, we prove some structural properties of the maximizer of $\theta$ with respect to its decomposition, and provide an extension of the approximation algorithm of Feige et al [4], which may be of independent interest.

## §2. Background

With V a finite set, let $\theta : 2^V \to \mathbb{R}$ be a set function which is a priori strictly bounded (in absolute value) by M, and (for technical simplicity) has a unique maximizer $V^* \in 2^V$. Note that this assumption is without loss of generality as a lexicographic ordering can be imposed upon $V' \in 2^V$ to distinguish between sets with equal function value. Alternately, one could imagine that by sacrificing an arbitrarily small additive difference from optimality, $\varepsilon > 0$, $\theta$ can be perturbed so that no two distinct sets have the same function value. We wish to maximize $\theta(V')$.

Recall that for a simple, undirected graph G = (V, E), the graph cut function $\delta: 2^V \to \mathbb{R}$ is defined as $\delta(V') = \left|\{(v_i, v_j) \in E \mid v_i \in V', v_j \in V \backslash V'\}\right|$. Similarly, the bivariate cut function $\delta_b : 2^V \times 2^V \to \mathbb{R}$ is defined as
$\delta_b(A, B) = \left|\{(v_i, v_j) \in E \mid v_i \in A, v_j \in B\}\right| \; \forall \, A, B \in 2^V$. We often make use of the following elementary fact.

**Fact 1:** Let G = (V, E) with $A, B \in 2^V$.
Then $\delta(A) + \delta(B) = \delta(A \cup B) + \delta(A \cap B) + 2\delta_b(A \backslash B, B \backslash A)$.

Proof: Let d(v) denote the degree of $v \in V$. For any $V' \in 2^V$ let $\omega(V')$ denote the size of the induced edge set |E(V')|, i.e. $|\{(u,v) \in E \mid u \text{ and } v \in V'\}|$. Observe that
$\delta(V') = \sum_{v \in V'} d(v) - 2\omega(V')$ (i). Hence for any $A, B \in 2^V$ we have

$\delta(A) + \delta(B) = \sum_{v \in A} d(v) + \sum_{v \in B} d(v) - 2\omega(A) - 2\omega(B)$

$= \sum_{v \in A \cap B} d(v) + \sum_{v \in A \cup B} d(v) - 2[\omega(A) + \omega(B)]$ (ii). Furthermore,

$\omega(A) + \omega(B) = 2\omega(A \cap B) + \delta_b(A \backslash B, A \cap B) + \delta_b(B \backslash A, A \cap B)$. Noting that
$\omega(A \cup B) = \omega(A \cap B) + \delta_b(A \backslash B, A \cap B) + \delta_b(B \backslash A, A \cap B) + \delta_b(A \backslash B, B \backslash A)$, we have
$\omega(A) + \omega(B) = \omega(A \cap B) + \omega(A \cap B) - \delta_b(A \backslash B, B \backslash A)$ (iii). Plugging (iii) into (ii), and using formula (i) to simplify, we get $\delta(A) + \delta(B) = \delta(A \cup B) + \delta(A \cap B) + 2\delta_b(A \backslash B, B \backslash A)$, as desired.

Recall that, by definition, $f: 2^V \to \mathbb{R}$ is *submodular* if
$\forall A, B \in 2^V \; f(A) + f(B) \geq f(A \cup B) + f(A \cap B)$. A function f is *supermodular* if (–f) is submodular, and *modular* if it is both sub- and super- modular. (Observe that since $\delta_b$ is non-negative, Fact 1 implies the well-known result that $\delta$ is submodular.) The following fact is also elementary.

**Fact 2:** Let $f_1, f_2, f_3 : 2^V \to \mathbb{R}$ be submodular, supermodular, and modular, respectively. Then $f_1 \pm f_3$ is submodular, $f_2 \pm f_3$ is supermodular, $f_1 - f_2$ is submodular, and $f_2 - f_1$ is supermodular.

A "decomposition" property that is essential for our purposes is that $\theta$ can be decomposed as $\alpha * (f - \delta)$, where f is submodular, $\delta$ is a cut function, and $\alpha$ is a positive constant.

**Prop. 3:** Let $\theta: 2^V \to \mathbb{R}$, then for some $\alpha > 0$, $\frac{1}{\alpha}\theta = f - \delta$, where f is submodular and $\delta$ is a cut function.

Proof: Let $M_1 = \max_{A,B \in 2^V}[\theta(A \cup B) + \theta(A \cap B) - \theta(A) - \theta(B)]/\alpha$. For sufficiently large positive $\alpha$, $M_1 \leq 1$. Let G be the complete graph $K_n$. Then $\delta_b(A\backslash B, B\backslash A) \geq 1$ whenever A and B are incomparable ($A \not\subset B$ and $B \not\subset A$), and therefore, by Fact 1, $\delta$ is not only submodular but in fact *strictly submodular* in the sense that $\delta(A) + \delta(B) > \delta(A \cup B) + \delta(A \cap B)$ whenever A and B are incomparable. It now follows from Fact 1 and the scaling $M_1 \leq 1$ that $f = \frac{1}{\alpha}\theta + \delta$ is also (strictly) submodular. □

In general, finding a decomposition of $\theta$ is nontrivial as we do not know $\alpha$. However, if $\theta$ is bounded by (the known constant) M, then clearly $\alpha = 4M$ will suffice. But rarely will we need $\alpha$ so large, as we may often decompose $\frac{1}{\alpha}\theta$ into $f - \delta$, even if $\frac{M_1}{\alpha} \gg 1$ [9]. For some $\theta$, $\frac{1}{\alpha}\theta$ is decomposable for an $\alpha \in (0,1)$. For simplicity, throughout the rest of this paper we assume that $\theta$ has already been rescaled, and so we shorten $\frac{1}{\alpha}\theta$ to $\theta$.

Next observe that more than one graph G may induce a cut function $\delta$ suitable for use in a decomposition of $\theta$. Indeed, a decomposition of $\theta$ into $f - \delta$, where $\delta$ is induced by a graph that is relatively sparse (and therefore has relative few dominating sets) is preferable due to the following structural property of the maximizer of $\theta$.

**Prop. 4:** Let V* be the unique maximizer of $\theta$. Then one of the following is true:
i) V* = $\varnothing$, ii) V* is a singleton, iii) V* is the unique maximizer of a function $\tilde{\theta}$, easily evaluable in terms of $\theta$, such that in any decomposition $f - \delta$ of $\tilde{\theta}$, V* is a dominating set in G.

Proof: First define $\tilde{\theta}: 2^V \to \mathbb{R}$ by $\tilde{\theta}(\varnothing) = \max\{\theta(\varnothing), \max_{v \in V}\theta(v)\}$, $\tilde{\theta}(V') = \theta(V') \ \forall \ V' \in 2^V \backslash \varnothing$. Assuming the presence of a value-giving oracle for $\theta$, we can evaluate $\tilde{\theta}$ in O(n) time. Clearly, every maximizer of $\theta$ is a maximizer of $\tilde{\theta}$, and $\tilde{\theta}$ has a unique maximizer unless $\theta$ is maximized on a singleton. Now suppose that none of i), ii), or iii) hold; then V* is a non-empty set and is non-dominating in some graph G decomposing $\tilde{\theta}$. Since V* is non-dominating, $\exists v \in V \backslash V^*$ such that $\delta_b(V^*, v) = 0$. As f is submodular, we have $f(V^* \backslash v) + f(v) \geq f(V^*) + f(\varnothing)$, which is equivalent to $\tilde{\theta}(V^* \backslash v) + \delta(V^* \backslash v) + \tilde{\theta}(v) + \delta(v) \geq \tilde{\theta}(V^*) + \delta(V^*) + \tilde{\theta}(\varnothing) + \delta(\varnothing)$. Applying Fact 1, we have $0 = 2\delta_b(V^* \backslash v, v) \geq [\tilde{\theta}(V^*) - \tilde{\theta}(V^* \backslash v)] + [\tilde{\theta}(\varnothing) - \tilde{\theta}(v)]$. But this is a contradiction since in the right hand side of the last inequality, the first bracketed term is positive because the maximizer of $\tilde{\theta}$ is unique, while the second bracketed term is nonnegative by the definition of $\tilde{\theta}(\varnothing)$. □

Prop. 4 asserts that if a sparse decomposition for $\tilde{\theta}$ were *given*, then the maximization of $\theta$ could be greatly simplified. Simply find a set $\tilde{V}$ maximizing $\tilde{\theta}$ over the dominating sets of G, and then take V* as $\text{argmax}\{\theta(V') \mid V' = \varnothing, \text{ a singleton, or } \tilde{V}\}$. It is not always easy to find a good (sparse) decomposition of $\theta$. In fact, it can be shown that any instance of the set covering problem corresponds to finding a sparsest graph G to decompose $\theta$, for some appropriate choice of $\theta$ [9]. At any rate, from this point on we assume that we are maximizing an already decomposed function $f - \delta$.

Additionally, Prop. 4 does not specify *how* to maximize $\tilde{\theta}$ over the dominating sets of G. Consequently, it is seemingly impractical to maximize $\theta$ unless $\theta$ is decomposed in such a way that the dominating sets of the attendant graph G are few and/or highly structured. In the next section we propose an alternate optimization procedure that is more clearly detailed and also less directly dependent upon the choice of G. The price for increased specification and decreased structural dependence is paid in solving many submodular maximization problems, rather than performing a single optimization.

Maximizing a submodular function is non-trivial, but can be treated by the dichotomy algorithm of Goldengorin [7]. In §4, we explore conditions mitigating the difficulty of solving such problems. For now, we detail Goldengorin's method. The dichotomy algorithm is a branch and bound method, which uses the submodular preservation rules (stated in simplest form in Prop. 5) as a pruning criterion. Some notation is in order: for $\varnothing \subset V_1 \subset V_2 \subset V$, let interval $[V_1, V_2] = \{V' \in 2^V \mid V_1 \subset V' \subset V_2\}$. Furthermore, for any $g : 2^V \to \mathbb{R}$, let $g^*[V_1, V_2]$ denote $\max_{V' \in [V_1, V_2]} g(V')$.

**Prop. 5 (Preservation Rules):** Let g be a submodular function on the interval $[V_1, V_2] \subset [\varnothing, V]$, and let $v \in V_2 \setminus V_1$. Then the following assertions hold.
(a) $g^*[V_1, V_2 \setminus v] - g^*[V_1 \cup v, V_2] \geq g(V_1) - g(V_1 \cup v)$.
(b) $g^*[V_1 \cup v, V_2] - g^*[V_1, V_2 \setminus v] \geq g(V_2) - g(V_2 \setminus v)$.

Proof: See [7].

The preservation rules translate into a branch and bound algorithm in the following straightforward manner. Namely, given an interval $[V_1, V_2]$ of $[\varnothing, V]$ over which we wish to maximize g, consider the $2|V_2 \setminus V_1|$ branches corresponding to subintervals of the form $[V_1, V_2 \setminus v]$ and $[V_1 \cup v, V_2]$. For any $v \in V_2 \setminus V_1$, if the right-hand side $g(V_1) - g(V_1 \cup v)$ of (a) Prop. 5 is non-negative then we may prune the interval $[V_1 \cup v, V_2]$, and/or if the right hand side $g(V_2) - g(V_2 \setminus v)$ of (b) is non-negative then we may prune the interval $[V_1, V_2 \setminus v]$.

Unfortunately, it is not simple to extend these preservation rules to maximizing $f - \delta$ over $[V_1, V_2]$, because f and $\delta$ are unlikely to be maximized/minimized on the same member of the pair $\{[V_1, V_2\backslash v], [V_1 \cup v, V_2]\}$ of intervals. The primary contribution of this paper is to show how to maximize $f - \delta$ by solving a sequence of submodular maximization problems. Because the dichotomy algorithm can be used to maximize *any* submodular function (as opposed to common restrictions such as monotonicity or non-negativity), it can be exploited as a subroutine in our branch and bound method to ensure the latter's convergence to a global maximizer of $f - \delta$. Solving an intractable problem with each call to a subroutine sounds daunting. However, as aforementioned, in §4 we show how to simplify the solution of these submodular maximization problems to greatly reduce the algorithm's overall complexity.

### §3. A Branch and Bound Maximization Procedure

The branch and bound (BB) algorithm takes as input a submodular function f, a cut function $\delta$ (with corresponding graph G), and outputs a maximizer of $f - \delta$. For simplicity, we assume that $f - \delta$ has a unique maximizer V*, and that f is normalized, i.e. $f(\emptyset)=0$. The latter assumption is without loss of generality as subtracting the constant $f(\emptyset)$ from f (and hence $f - \delta$) will yield a normalized function. A rooted tree is then constructed, where each node of the tree has a different attendant graph $\hat{G} = \left( V, \hat{E} \subset \binom{V}{2} \right)$, with the graph corresponding to the root node being $\hat{G} = (V, \emptyset)$. (Recall that $\binom{V}{2}$ denotes $\{(v_i, v_j) \mid v_i, v_j \in V, v_i \neq v_j\}$.) Choosing $f_u$ as a submodular upper bound of f, agreeing with f on $\emptyset$ and singletons, at each node of the tree we solve the relaxed subproblem $P(\hat{G}) = \{\max (f_u - \hat{\delta})(V') \mid V' \in \hat{F}\}$, where $\hat{\delta}(V') = \sum_{v \in V'} \hat{d}(v)$ with $\hat{d}(v) = |\{(u,v) \in E \cap \hat{E}\}|$, and $\hat{F}$ is a subfamily of the independent sets of $\hat{G}$.

Solving $P(\hat{G})$ yields $\hat{V}_1$, a maximizer of $f_u - \hat{\delta}$ over $\hat{F}$, and $\hat{V}_2$ a maximal extension of $\hat{V}_1$ in $\hat{F}$. It will be shown (Prop. 9) that $f_u - \hat{\delta}$ is an upper bound of $f - \delta$ on $\hat{F}$, hence $(f_u - \hat{\delta})(\hat{V}_1) \geq (f-\delta)(V') \; \forall V' \in \hat{F}$. Denoting the $V' \in 2^V$ with the largest $f - \delta$ value thus far observed as Z, we "prune" the current node if $(f_u - \hat{\delta})(\hat{V}_1) \leq (f-\delta)(Z)$, and "fathom" it if $\max_{i=1,2}(f-\delta)(\hat{V}_i) = (f_u - \hat{\delta})(\hat{V}_1) > (f-\delta)(Z)$. Otherwise, we grow the tree by adding child nodes to the current node, each of which corresponds to a previously unseen graph $\hat{G}'$. At every node, the $\hat{G}'$ are chosen from the same finite family and possess significant structure (see Fact 6). As a result, the BB algorithm must terminate, and the specified structure of the $\hat{G}$ (see Prop. 10) guarantee that the final Z output is V*.

Our program variables, $\hat{E}$, Z, and z* denote the current subset of $\binom{V}{2}$ under consideration, the best observed maximizer of $f - \delta$, and $(f - \delta)(Z)$, respectively.

**Alg. BB**

Input: A function f-$\delta$ to be maximized, where f is normalized ( f($\emptyset$)=0 ) and submodular, and $\delta$ is the cut function for a graph G = (V, E).

Output: A maximizer V* of f-$\delta$.

Initialization: Z = $\emptyset$, $\hat{E}$ = $\emptyset$, z* = f($\emptyset$) - $\delta$($\emptyset$). Label the initial node, corresponding to
$\hat{G}$ = (V, $\emptyset$), open.

Procedure:

1) At the open node corresponding to the graph $\hat{G} = \left( V, \hat{E} \subset \binom{V}{2} \right)$, solve the "relaxed" subproblem $P(\hat{G}) = \{\max (f_u - \hat{\delta})(V') \mid V' \in \hat{F}\}$. $\hat{F}$, denoting the feasible region for $P(\hat{G})$, is defined as $2^{I(\hat{G})}$, where $I(\hat{G})$ is the unique maximum size independent set of $\hat{G}$, if $\hat{G} \neq K_n$, and $\hat{F} = \{\emptyset, v_1, ..., v_n\}$ if $\hat{G} = K_n$. In both cases, $f_u(V') \geq f(V')$ with equality holding if V' is a singleton or $\emptyset$, and $\hat{\delta}(V') = \sum_{v \in V'} \hat{d}(v)$ with $\hat{d}(v) = |\{(u,v) \in E \cap \hat{E}\}|$.

Solve $P(\hat{G})$ to obtain the solution $\hat{V}_1$. *(Relaxing)*

2) Set $\hat{V}_2 = I(\hat{G})$ if $\hat{G} \neq K_n$, else set $\hat{V}_2 = \hat{V}_1$.

3) Set $\hat{\theta}_1 = (f_u - \hat{\delta})(\hat{V}_1)$ and $\hat{\theta}_2 = \max_{i=1,2} (f - \delta)(\hat{V}_i)$. *(Bounding)*

4)     a) If $\hat{\theta}_1 \leq z^*$, then label the node corresponding to $\hat{G}$ closed. *(Pruning)*

       b) If $\hat{\theta}_2 > z^*$, update Z $\to$ argmax$_{i=1,2}$ $(f-\delta)(\hat{V}_i)$ and set z* $\to \hat{\theta}_2$. Furthermore, if $\hat{\theta}_1 = \hat{\theta}_2$, label the node corresponding to $\hat{G}$ closed. *(Fathoming)*

       c) Else ( $\hat{\theta}_2 \leq z^* < \hat{\theta}_1$ ), consider $|I(\hat{G})|$ new nodes corresponding to graphs of the form $\hat{G} \setminus v_i = \left( V, \hat{E} \cup \{(v_i, v_j) \mid v_j \in V, v_j \neq v_i\} \right)$, where $v_i \in I(\hat{G})$. For each graph $\hat{G} \setminus v_i$, check if $\hat{G} \setminus v_i$ has occurred as the attendant graph to some node previously created by BB, and if not, create a new node with graph $\hat{G} \setminus v_i$ and label that node open. Label the current node, corresponding to $\hat{G}$, closed. Select any open node and go to step 1.
*(Branching)*

5) When all nodes have been closed, output Z.

      We call the graphs $\hat{G}$ that arise during BB *astrals*, as (Fact 6 below) each $\hat{G} \neq K_n$ is the union of stars centered at $v \in V \setminus I(\hat{G})$. Throughout the rest of this paper, we refer explicitly and implicitly to the following structural properties of astrals.

**Fact 6**: Let $\hat{G} \neq K_n$ be an astral. Then:

i) $\hat{G}$ has a unique maximal independent set of size $\geq 2$.

ii) Any *maximal* independent set of $\hat{G}$ with size $\geq 2$ is also maximum.
iii) $\hat{G}$ is the union of stars centered at $v \in V \setminus I(\hat{G})$.

Each claim in Fact 6 can be verified by induction. For convenience, we refer to an independent set of size 0 or 1 as "trivial" as such sets are always independent in loopless graphs.

The following two observations are of importance in establishing the correctness of BB and bounding the number of nodes visited. We say that a graph $\hat{G}$ is reachable by BB if, in the absence of pruning or fathoming, BB would have visited a node for which $\hat{G}$ is the attendant graph.

**Prop. 7**: Let $V' \in 2^V$ such that $|V'| \geq 2$. Then there exists $\hat{G}$ reachable by BB for which $I(\hat{G}) = V'$.

Proof: Suppose not, and let V' be a maximum size element of $2^V$ for which there is no $\hat{G}$ reachable by BB with $I(\hat{G}) = V'$. Since the graph corresponding to the root node is $\hat{G} = (V, \emptyset)$, we have $V' \subset_{\neq} V$. Let $v_i \in V \setminus V'$ and note that by assumption, $V' \cup v_i = I(\tilde{G})$ for some $\tilde{G}$ reachable by BB. Consider the graph $\tilde{G} \setminus v_i = (V, \tilde{E} \cup \{(v_i, v_j) \mid v_j \in V, v_j \neq v_i\})$. In the absence of pruning or fathoming, $\tilde{G} \setminus v_i$ must arise as the graph attendant to a child node of the graph corresponding to $\tilde{G}$, or as the graph attendant to some node previously visited by BB. In either case, $\tilde{G} \setminus v_i$ is reachable by BB.

Because $|V'| \geq 2$, $|V' \cup v_i| = |I(\tilde{G})| \geq 3$. Thus $I(\tilde{G}) \setminus v_i = V'$ is an independent set of $\tilde{G} \setminus v_i$ of size $\geq 2$. Furthermore, since (Fact 6 i)) $I(\tilde{G})$ is the only non-trivial maximal independent set of $\tilde{G}$, it follows that any non-trivial maximal independent set of $\tilde{G} \setminus v_i$, in particular $I(\tilde{G} \setminus v_i)$, must be a subset of $I(\tilde{G}) \setminus v_i$. Thus $I(\tilde{G} \setminus v_i) = I(\tilde{G}) \setminus v_i = V'$. Hence $V' = I(\hat{G})$ for some $\hat{G}$ reachable by BB, contradicting our original assumption.

**Prop. 8**: $I : \{\hat{G} \mid \hat{G} \text{ is an astral} \neq K_n\} \to 2^V$ is injective.
Proof: Let $\hat{G}_1$ and $\hat{G}_2$ be astrals, neither of which = $K_n$, such that $I(\hat{G}_1) = I(\hat{G}_2)$. Then, by Fact 6 iii), this implies
$\hat{G}_1 = (V, \bigcup_{v_i \notin I(\hat{G}_1)} \{(v_i, v_j) \mid v_j \in V, v_j \neq v_i\}) = (V, \bigcup_{v_i \notin I(\hat{G}_2)} \{(v_i, v_j) \mid v_j \in V, v_j \neq v_i\}) = \hat{G}_2$.

Step 4c guarantees that each node created corresponds to a different astral, while Prop. 7 and Prop. 8 imply that every astral $\hat{G} \neq K_n$ corresponds to a different, non-empty, element of $2^V$. Hence, at most $2^n - n$ nodes will be visited by BB. We may extend I to an injection on all astrals by defining $I(K_n) = \emptyset$. Since each $\hat{G}$ now corresponds to a different member of $2^V$, we may perform the "previously created" evaluation in step 4c in unit time (but exponential space) by creating a table with an entry for each $V' \in 2^V$, and marking that entry if V' has occurred as $I(\hat{G})$ for some $\hat{G}$ created by BB. Note that our space requirements are naïve and can be greatly reduced by choosing an appropriate branching scheme. For ease of notation, we shall now refer to a node and its attendant graph interchangeably.

Since BB visits a finite number of nodes, it must converge. Letting $Z_F$ denote the Z output by BB, we need to show that $Z_F = V^*$.

**Prop. 9**: $f_u - \hat{\delta}$ is an upper bound of $f - \delta$ on $\hat{F}$.

Proof: Let $V' \in \hat{F}$. By definition of $f_u$, we have $f_u(V') \geq f(V')$. Thus it is sufficient to show that $\hat{\delta}(V') \leq \delta(V')$ to prove the claim. To see this, recall that
$\hat{\delta}(V') = \sum_{v \in V'} \hat{d}(v) = \sum_{v \in V'} |\{(u,v) \in \hat{E} \cap E\}|$. Since $V' \in \hat{F}$, V' is independent in $\hat{G}$ and so
$\forall v \in V'$ $\{(u,v) \in \hat{E} \cap E\} = \{(u,v) \in \hat{E} \cap E \mid u \notin V'\} \subset \{(u,v) \in E \mid u \notin V'\}$. Thus
$\hat{\delta}(V') \leq \sum_{v \in V'} |\{(u,v) \in E \mid u \notin V'\}|$. Now, for any $v_1, v_2 \in V'$ with $v_1 \neq v_2$ we have
$\{(u, v_1) \in E \mid u \notin V'\} \cap \{(u,v_2) \in E \mid u \notin V'\} = \emptyset$, giving us
$\sum_{v \in V'} |\{(u,v) \in E \mid u \notin V'\}| = |\bigcup_{v \in V'} \{(u,v) \in E \mid u \notin V'\}|$. The right hand side of the preceding equality is equal to $|\{(u,v) \in E \mid u \notin V', v \in V'\}| = \delta(V')$. Thus $\hat{\delta}(V') \leq \delta(V')$, as desired.
□

**Prop. 10**: $Z_F = V^*$.

Proof: If $V^*$ ever arose as $\hat{V}_1$ or $\hat{V}_2$ for some $P(\hat{G})$ reached during BB, or as our initial $Z = \emptyset$, then $Z_F = V^*$. Otherwise we must have $f(Z_F) - \delta(Z_F) > f(V^*) - \delta(V^*)$, contradicting the fact that $V^*$ is a maximizer of $f - \delta$. Now, $V^*$ is either: $\emptyset$, a singleton, or has size $\geq 2$. In the first case, $Z_F = V^*$ as we initialize $Z = \emptyset$. For $V^*$ of size $> 0$, let $G^*$ denote $K_n$, if $V^*$ is a singleton, or $G^* = (V, E^*)$ where
$E^* = \left\{ (v_i, v_j) \in \binom{V}{2} \mid \text{at least one of } v_i, v_j \notin V^* \right\}$, if $|V^*| \geq 2$. First we show that if we reach $P(G^*)$, then we must have $Z_F = V^*$.

If $V^*$ is a singleton, then $V_1^* \in \{\emptyset, v_1, ..., v_n\}$. Since $f_u = f$ on this set, and $\delta^* = \delta$ on $\{\emptyset, v_1, ..., v_n\}$ for $G^* = K_n$, we see that $f_u - \delta^* = f - \delta$ on this set. Thus $V_2^* = V_1^* = V^*$. If $|V^*| \geq 2$, then $V^* = I(G^*)$, and so $V_2^* = V^*$. In both cases, $\theta_2^* \geq z^*$. Either $\theta_2^* > z^*$, in which case we update $Z \to V^*$, or $\theta_2^* = z^*$, in which case Z was already equal to $V^*$ since $f - \delta$ has a unique maximizer. In either case, Z, and thus $Z_F$, is equal to $V^*$.

Notice that we only fail to reach $P(G^*)$ if some previous graph $\hat{G} = (V, \hat{E} \subset E^*)$ for which $V^* \in \hat{F}$ was fathomed or pruned. This can happen in one of two ways, either $\hat{\theta}_1 = \hat{\theta}_2 >$ the (then) value of $z^*$ (fathoming), or $\hat{\theta}_1 \leq$ the (then) value of $z^*$ (pruning). Observe that $(f_u - \hat{\delta})(V^*) \geq (f_u - \delta^*)(V^*)$. This follows since $\hat{E} \subset E^*$ and so $\hat{\delta}(V^*) \leq \delta^*(V^*)$. Thus $\hat{\theta}_1 \geq (f_u - \hat{\delta})(V^*) \geq (f_u - \delta^*)(V^*) \geq (f - \delta)(V^*)$, as $f_u - \delta^* \geq f - \delta$ on $F^*$. In the first case, by the case assumption, we have $\hat{\theta}_2 = \hat{\theta}_1$ and we must have updated $Z \to \text{argmax}_{i=1,2} (f-\delta)(\hat{V}_i)$, in which case $f(Z) - \delta(Z) = \hat{\theta}_2 = \hat{\theta}_1 \geq f(V^*) - \delta(V^*)$, and so $Z$ must equal $V^*$.

In case two, if $z^* = f(V^*) - \delta(V^*)$ then $Z$ must already be set equal to $V^*$. Otherwise, $\hat{\theta}_1 \leq z^* < f(V^*) - \delta(V^*)$, but this contradicts the previous inequality. Thus in every scenario where $P(G^*)$ is not reached that does not lead to a contradiction, we must have set $Z = V^*$. Hence in all cases (whether $P(G^*)$ is reached or not) BB must output $Z_F = V^*$. □

### §4. Improvements and Properties of the Procedure

Because $f_u$ is submodular and $\hat{\delta}$ is modular, by Fact 2, solving each subproblem $P(\hat{G})$ of BB involves the maximization of a submodular function. While such problems are, in general, intractable, they possess significant structure and can be solved by the dichotomy algorithm. Since the size of the ground set, $I(\hat{G})$, shrinks as we increase distance from the root node, using the dichotomy algorithm to solve $P(\hat{G})$ may indeed be practical for nodes sufficiently distant from the root. In this section we discuss how to exactly or approximately solve $P(\hat{G})$ via judicious choice of $f_u$, and approximation methods, respectively. These techniques can be used to treat subproblems $P(\hat{G})$ where $2^{I(\hat{G})}$ is too large for a practical implementation of the dichotomy algorithm. We also compare the advantages and disadvantages of treating $P(\hat{G})$ by these two methods (choice of $f_u$ and approximation). The section is concluded by a proof that BB is interruptible; that is, if BB is stopped before all nodes have been closed, one can bound the additive difference from optimality of the current best solution, Z.

Consider choosing $f_u$ as $f_u(V') = \sum_{v \in V'} f(v)$. Cleary $f_u(V') = f(V')$ if V' is a singleton or $\varnothing$ (as $f(\varnothing)=0$). Because f is submodular, we have $f_u(V') \leq \sum_{v \in V'} f(v)$, and so $f_u$ is an upper bound of f as well. Finally, $f_u$ is clearly modular and thus submodular, fulfilling all of the requirements of BB.

**Prop. 11**: For $f_u$ defined as $f_u(V') = \sum_{v \in V'} f(v)$, solving $P(\hat{G})$ requires $O(n)$ effort.

Proof: For $f_u(V') = \sum_{v \in V'} f(v)$, $f_u - \hat{\delta}$ is a modular function. If $\hat{G} = K_n$, the claim is independent of the particular choice of $f_u$. For $\hat{G} \neq K_n$, observe that $f_u - \hat{\delta}$ is maximized (over $2^{I(\hat{G})}$) on $\hat{V} = \cup \{v \in I(\hat{G}) \mid (f_u - \hat{\delta})(v) \geq 0\}$. Clearly, finding $\hat{V}$ requires $O(n)$ effort.

Choosing $f_u(V') = \sum_{v \in V'} f(v)$ has the advantage of bounding the complexity of BB by $n2^n$. However, it may happen that $f_u(V') \ll \sum_{v \in V'} f(v)$ for most $V' \in 2^V$. In which case, pruning and fathoming will be rare. Now consider selecting $f$ as $f_u$. Because $f$ is an arbitrary submodular function, solving $P(\hat{G})$ will be NP-hard. However, $f$ is clearly the tightest upper bound to itself, and furthermore, $f - \hat{\delta}$ will converge to $f - \delta$ when $\hat{G}$ is sufficiently dense. Specifically, $(f - \hat{\delta})(V') = (f - \delta)(V')$ if $V'$ is independent in $\hat{G}$ (is a subset of $I(\hat{G})$) and all of the cut edges of $V'$ (in $G$) are present in $\hat{E}$. Since $f - \hat{\delta}$ converges monotonically to $f - \delta$, pruning will also be more frequent. We attempt to maintain the advantages of a tighter choice of $f_u$, while mitigating the difficulty of solving an NP-hard problem, by introducing the LS algorithm of Feige et al. [4] to approximately solve $P(\hat{G})$. BB can then be adapted to use these approximate solutions while still converging to $V^*$.

Before proceeding to the LS algorithm, note that we may use either choice of $f_u$ when running BB and we will still converge to $V^*$, even alternating $f_u$ throughout the course of the algorithm. For example, we may take $f_u(V') = \sum_{v \in V'} f(v)$ for $P(\hat{G})$ if $f$ is approximately modular on $2^{I(\hat{G})}$, and occasionally use $f_u = f$ in hope of generating larger $\hat{\theta}_2$ values, resulting in increased pruning.

**Alg. LS**
Input: A non-negative submodular function $g: 2^X \to \mathbb{R}$, $\varepsilon > 0$.
Output: A $\left(\frac{1}{3} - \frac{\varepsilon}{r}\right)$-optimal maximizer of $g$, where $r = |X|$.

Procedure:
1) Let $X' = \text{argmax}_{x \in X} g(x)$.
2) If $\exists x \in X \setminus X'$ such that $g(X' \cup x) > \left(1 + \frac{\varepsilon}{r^2}\right) g(X')$, update $X' \to X' \cup x$ and go back to step 2.
3) If $\exists x \in X'$ such that $g(X' \setminus x) > \left(1 + \frac{\varepsilon}{r^2}\right) g(X')$, update $X' \to X' \setminus x$ and go back to step 2.
4) Output max $g(X')$, $g(X \setminus X')$.

**Prop. 12**: Let $g: 2^X \to \mathbb{R}$ be a non-negative submodular function. Then LS returns a $\left(\frac{1}{3} - \frac{\varepsilon}{r}\right)$-optimal maximizer of g in $O\left(\frac{1}{\varepsilon} r^3 \log r\right)$ operations.

Proof: See [4].

Now, observe that if $f - \delta$ is non-negative, then so is $f_u - \hat{\delta}$ when restricted to $2^{I(\hat{G})}$, as $f_u - \hat{\delta} \geq f - \delta$ on this domain by Prop. 9. Hence, taking $\varepsilon = 1$ and assuming $|I(\hat{G})| \geq 12$, LS outputs a ¼-optimal solution of $P(\hat{G})$ in no greater than $O(n^3 \log n)$ time. BB can be adapted to use these approximate solutions, and still converge to V*, as follows. Replace the pruning condition $\hat{\theta}_1 \leq z^*$ with $4\hat{\theta}_1 \leq z^*$, and replace the fathoming condition $\hat{\theta}_1 = \hat{\theta}_2$ with $4\hat{\theta}_1 \leq \hat{\theta}_2$. Making the adjustments yields a modified BB with a greater (by a factor of $n^2 \log n$) worst case complexity, but which is more likely to have pruned nodes. LS can be modified to produce good heuristic solutions for submodular functions $g: 2^X \to \mathbb{R}$ that are not necessarily non-negative everywhere, but which are non-negative on a certain family of subsets of $2^X$. We take up this issue in the next section on constrained optimization. We conclude this section by establishing that BB is interruptible.

We say that a subproblem $P(\hat{G})$, or graph $\hat{G}$, has depth k if its attendant node has (directed) distance k from the root node. For simplicity of statement, we suppose that each subproblem of depth $\leq$ k has been closed.

**Prop. 13**: Suppose that each subproblem $P(\hat{G})$ of depth $\leq$ k has been closed (by branching, pruning, or fathoming). For each $\hat{G}$ of depth k reached by BB, let $\hat{\Delta} = \hat{\theta}_1 - \hat{\theta}_2$. Then $(f-\delta)(V^*) \leq (f-\delta)(Z) + \Delta_{max}$, where Z is the current best solution and $\Delta_{max} = \max \hat{\Delta}$.

Proof: Let G* be defined as in the proof of Prop. 10. By Prop. 10, after closing all subproblems of depth $\leq$ k, either $Z = V^*$, or $\hat{E} \subset E^*$ for some $\hat{G}=(V, \hat{E})$ of depth k. If $Z = V^*$, the claim immediately holds as $\hat{\theta}_1$ is always $\geq \hat{\theta}_2$. Else observe that $\hat{\theta}_1 \geq (f-\delta)(V^*)$ (by the proof of Prop. 10), and so $\hat{\theta}_2 + \Delta_{max} \geq \hat{\theta}_2 + (\hat{\theta}_1 - \hat{\theta}_2) \geq (f-\delta)(V^*)$. Since $\hat{\theta}_2 \leq (f-\delta)(Z)$, the desired result follows.

It is easy to see that this result holds without the assumption that all subproblems are of the same level. It is also clear that a similar result holds for when the $P(\hat{G})$ are solved approximately instead of exactly.

**§5. Extension to Constrained Problems**

In this section we extend our previous results to the constrained maximization of $\theta$ over a family $S \subset 2^V$. In particular, we consider the case where S is a *subset system*, a notion generalizing such well-known constraints as independence in a graph, and a collection of sets forming a packing. We show that problems so-constrained can be solved via reformulation into an equivalent unconstrained maximization problem, or by modifying algorithm BB. The latter approach is of interest if $\theta$ possesses desirable properties for unconstrained maximization, such as a sparse decomposition and non-negativity.

Given a family $\mathscr{F} \subset 2^V$, we say that $p: 2^V \to \mathbb{R}$ is a *membership function* of $\mathscr{F}$ if $V' \in \mathscr{F}$ iff $p(V') \leq 0$. A family $\mathscr{F} \subset 2^V$ is a subset system if $\forall\, V_1 \in \mathscr{F}$ and any $V_2 \subset V_1$, $V_2 \in \mathscr{F}$. In particular, the independent sets of a graph form a subset system. Subset systems and *super*modular functions are intimately related via the following proposition.

**Prop. 14:** Let $S \subset 2^V$ be a non-empty subset system, then there exists a non-negative, normalized ($q(\varnothing)=0$), supermodular function q such that q is a membership function of S. Conversely, if q is a non-negative, normalized, supermodular function, then $S = \{V' \in 2^V \mid q(V') \leq 0\}$ is a non-empty subset system.

Proof: Let S be a non-empty subset system, and consider $q: 2^V \to \mathbb{R}$ defined by $q(V') = \begin{cases} 0, & \text{if } V' \in S \\ \exp(|V'|), & \text{if } V' \notin S \end{cases}$. Clearly q is normalized, and $q(V') \leq 0$ iff $V' \in S$. To prove q is supermodular we must show that $\forall A, B \in 2^V\ q(A) + q(B) \leq q(A \cup B) + q(A \cap B)$. This is obviously true if $A \subset B$ or $B \subset A$, so let A and B be incomparable. If $A, B \in S$, $q(A) = q(B) = 0$, and the inequality follows trivially. If $A \in S, B \notin S$, then $A \cap B \in S$ and $A \cup B \notin S$, so that $q(A) + q(B) = \exp(|B|) \leq \exp(|A \cup B|) = q(A \cup B) + q(A \cap B)$. Similarly, if $A \notin S, B \in S$. Finally, if both $A \notin S, B \notin S$, then $q(A) + q(B) = \exp(|A|) + \exp(|B|) \leq 2\exp(\max\{|A|, |B|\}) \leq \exp(\max\{|A|, |B|\} + 1)$ (since $e > 2$) $\leq \exp(|A \cup B|)$ (since A and B are incomparable) $\leq q(A \cup B) + q(A \cap B)$ (since $A \cup B \notin S$). Thus q is supermodular.

Conversely, let q be non-negative, normalized, and supermodular. Let $S = \{V' \in 2^V \mid q(V') \leq 0\}$; since q is normalized, $q(\varnothing) = 0$, and so S is non-empty. Now let $V_1 \in S$ and $V_2 \subsetneq V_1$. Then we have $q(V_2) + q(V_1 \setminus V_2) \leq q(V_1) + q(\varnothing) = q(V_1)$, by the normalization and supermodularity of q. Since $V_1 \in S$, this inequality implies $q(V_2) \leq q(V_1) \leq 0$, and so $V_2 \in S$. $\square$

Next, an immediate consequence of Prop. 14 is that the maximization of $\theta$ over a non-empty subset system, S, is equivalent to the unconstrained maximization of a certain set function.

**Cor. 15:** Let $(P) = \{\max \theta(V') \mid V' \in S\}$, where $\theta$ is strictly bounded in absolute value by M, and S is a non-empty subset system. Then (P) is equivalent to the unconstrained maximization problem $(P') = \{\max \theta(V') - Mq(V') \mid V' \in 2^V\}$, where q is as defined in (the first part of the proof of) Prop. 14.

Proof: Let V* maximize $\theta$ over S, and notice that $(\theta - Mq) = \theta$ $\forall V' \in S$, by definition of q. And so, $(\theta - Mq)(V^*) \geq (\theta - Mq)(V')$ $\forall V' \in S$. Furthermore, for any $V' \in 2^V \setminus S$ we have $q \geq e$ and therefore, since $M > \max_{V' \in 2^V} |\theta(V')|$, $(\theta - Mq)(V') < \min_{V' \in 2^V} \theta(V')$. Thus $(\theta - Mq)(V^*) = \theta(V^*) > (\theta - Mq)(V')$ $\forall V' \in 2^V \setminus S$, and so V* is a global maximizer of $(\theta - Mq)$.
    Conversely, let $\tilde{V}$ maximize $(\theta - Mq)$ over $2^V$. Since $(\theta - Mq)$ agrees with $\theta$ over S, this implies $(\theta - Mq)(\tilde{V}) \geq (\theta - Mq)(V^*) = \theta(V^*)$. By our earlier observation, $\theta(V^*) > (\theta - Mq)(V')$ $\forall V' \in 2^V \setminus S$. Thus we must have that $\tilde{V} \in S$, which further implies that $\theta(\tilde{V}) = (\theta - Mq)(\tilde{V}) \geq \theta(V^*)$, and so $\tilde{V}$ is a maximizer of $\theta$ over S. □

Note that if our original function to be maximized over S, $\theta$, were submodular, then so is the equivalent unconstrained function $(\theta - Mq)$, by Fact 2.

    Clearly Cor. 15 provides a method for solving (subset system) constrained maximization problems, by converting them to equivalent unconstrained problems to which the results of the previous sections apply. Suppose, however, that $\theta$ possesses some properties that are useful in the unconstrained maximization case, such as a sparse decomposition and non-negativity. It is reasonable to assume that these properties would also be of use when maximizing $\theta$ over S, but they will be lost due to the reformulation. Hence we are motivated to adapt our first maximization algorithm BB to a form designed specifically for constrained problems, BBC. Similarly to before, we assume that $\theta$ is pre-decomposed as $f - \delta$, and possesses a unique maximizer, V* *over S*. We validate this alternate maximization method by showing that if $f - \delta$ is non-negative, then each of the subproblems $P(\hat{G})$ of BBC is treatable with a heuristic adapted from LS.

**Alg. BBC**
Input: A function $f - \delta$ to be maximized, where f is normalized ($f(\emptyset) = 0$) and submodular, and $\delta$ is the cut function for a graph $G = (V, E)$. A non-empty subset system, S, over which to maximize $f - \delta$.

Output: A maximizer V* of $f - \delta$ over S.

Initialization: $Z = \emptyset$, $\hat{E} = \emptyset$, $z^* = f(\emptyset) - \delta(\emptyset)$. Label the initial node, corresponding to
        $\hat{G} = (V, \emptyset)$, open.
Procedure:

1) At the open node corresponding to the graph $\hat{G} = \left(V, \hat{E} \subset \binom{V}{2}\right)$, solve the "relaxed" subproblem $P(\hat{G}) = \{\max (f_u - \hat{\delta})(V') \mid V' \in \hat{F}\}$. Here $\hat{F} = S \cap 2^{I(\hat{G})}$, if $\hat{G} \neq K_n$, and $\hat{F} = S \cap \{\emptyset, v_1, \ldots, v_n\}$ if $\hat{G} = K_n$. In both cases, $f_u(V') \geq f(V')$ with equality holding if $V'$ is a singleton or $\emptyset$, and $\hat{\delta}(V') = \sum_{v \in V'} \hat{d}(v)$ with $\hat{d}(v) = |\{(u,v) \in E \cap \hat{E}\}|$. Solve $P(\hat{G})$ to obtain the solution $\hat{V}_1$. *(Relaxing)*

2) Using a greedy method, extend $\hat{V}_1$ to a maximal cardinality element $\hat{V}_2$ of $\hat{F}$.

3) Set $\hat{\theta}_1 = (f_u - \hat{\delta})(\hat{V}_1)$ and $\hat{\theta}_2 = \max_{i=1,2} (f - \delta)(\hat{V}_i)$. *(Bounding)*

4)  a) If $\hat{\theta}_1 \leq z^*$, then label the node corresponding to $\hat{G}$ closed. *(Pruning)*

   b) If $\hat{\theta}_2 > z^*$, update $Z \to \operatorname{argmax}_{i=1,2} (f - \delta)(\hat{V}_i)$ and set $z^* \to \hat{\theta}_2$. Furthermore, if $\hat{\theta}_1 = \hat{\theta}_2$, label the node corresponding to $\hat{G}$ closed. *(Fathoming)*

   c) Else ($\hat{\theta}_2 \leq z^* < \hat{\theta}_1$), consider $|I(\hat{G})|$ new nodes corresponding to graphs of the form $\hat{G} \setminus v_i = \left(V, \hat{E} \cup \{(v_i, v_j) \mid v_j \in V, v_j \neq v_i\}\right)$, where $v_i \in I(\hat{G})$. For each graph $\hat{G} \setminus v_i$, check if $\hat{G} \setminus v_i$ has occurred as the attendant graph to some node previously created by BB, and if not, create a new node with graph $\hat{G} \setminus v_i$ and label that node open. Label the current node, corresponding to $\hat{G}$, closed. Select any open node and go to step 1. *(Branching)*

5) When all nodes have been closed, output Z.

We assume the existence of an oracle to determine membership in S in polynomial time (i.e. that $p(V')$ can be evaluated in polynomial time). Similar to before, our branch and bound tree will have $O(2^n)$ nodes, and each V' of size $\geq 2$ will correspond to $I(\hat{G})$ for some graph $\hat{G}$ whose attendant node is reachable by BBC. In particular, for each $V' \in S$ of size $\geq 2$, this implies that $V' = \hat{F}$ for the subproblem $P(\hat{G})$ associated with $\hat{G}$. In a proof largely identical to that of Prop. 10, we show that BBC converges to V*.

**Prop. 16**: $Z_F = V^*$.

Proof: First, note that $Z_F = \emptyset$ or $\hat{V}_2$ for a subproblem $P(\hat{G})$ reached by the algorithm, by definition of our updating rules. Hence $Z_F \in S$. Now, if V* ever arose as $\hat{V}_2$ for some subproblem $P(\hat{G})$ reached by BBC, or as our initial $Z = \emptyset$, then $Z_F = V^*$. Otherwise we must have $Z_F \in S$ and $f(Z_F) - \delta(Z_F) > f(V^*) - \delta(V^*)$, contradicting the fact that V* is a maximizer of $f - \delta$ over S. Now, V* is either: $\emptyset$, a singleton, or has size $\geq 2$. In the first case, $Z_F = V^*$, as we initialize $Z = \emptyset$. For V* of size > 0, let G* denote $K_n$ if V* is a singleton, or $G^* = (V, E^*)$ where

$E^* = \left\{ (v_i, v_j) \in \binom{V}{2} \mid \text{at least one of } v_i, v_j \notin V^* \right\}$, if $|V^*| \geq 2$. First we show that if we reach $P(G^*)$, then we must have $Z_F = V^*$.

If $V^*$ is a singleton, then $V_1^* \in \{\emptyset, v_1, ..., v_n\}$. Since $f_u = f$ on this set, and $\delta^* = \delta$ on $\{\emptyset, v_1, ..., v_n\}$ for $G^* = K_n$, we see that $f_u - \delta^* = f - \delta$ on this set. Thus $V_2^* = V_1^* = V^*$. Now suppose $|V^*| \geq 2$. Observe that $V^* = I(G^*)$, and since $V_1^* \in F^* = 2^{I(G^*)} \cap S$, it follows that $V_1^* \subset V^*$. Since $V^* = I(G^*)$ and $V^* \in S$, we see that $V_1^*$ has a unique maximal extension in $F^*$, namely $V^*$. Hence $V_2^* = V^*$. In both cases, $\theta_2^* \geq z^*$. Either $\theta_2^* > z^*$, in which case we update $Z \to V^*$, or $\theta_2^* = z^*$, in which case $Z$ was already equal to $V^*$ since $f - \delta$ has a unique maximizer over $S$. In either case, $Z$, and thus $Z_F$, is equal to $V^*$.

Notice that we only fail to reach $P(G^*)$ if some previous graph $\hat{G} = (V, \hat{E} \subset E^*)$ for which $V^* \in \hat{F}$ was fathomed or pruned. This can happen in one of two ways, either $\hat{\theta}_1 = \hat{\theta}_2 >$ the (then) value of $z^*$ (fathoming), or $\hat{\theta}_1 \leq$ the (then) value of $z^*$ (pruning). Observe that $(f_u - \hat{\delta})(V^*) \geq (f_u - \delta^*)(V^*)$. This follows since $\hat{E} \subset E^*$ and so $\hat{\delta}(V^*) \leq \delta^*(V^*)$. Thus $\hat{\theta}_1 \geq (f_u - \hat{\delta})(V^*) \geq (f_u - \delta^*)(V^*) \geq (f - \delta)(V^*)$, as $f_u - \delta^* \geq f - \delta$ on $F^*$. In the first case, by the case assumption, we have $\hat{\theta}_2 = \hat{\theta}_1$ and we must have updated $Z \to \operatorname{argmax}_{i=1,2} (f - \delta)(\hat{V}_i)$, in which case $f(Z) - \delta(Z) = \hat{\theta}_2 = \hat{\theta}_1 \geq f(V^*) - \delta(V^*)$, and so $Z$ must equal $V^*$.

In case two, if $z^* = f(V^*) - \delta(V^*)$ then $Z$ must already be set equal to $V^*$. Otherwise, $\hat{\theta}_1 \leq z^* < f(V^*) - \delta(V^*)$, but this contradicts the previous inequality. Thus in every scenario where $P(G^*)$ is not reached that does not lead to a contradiction, we must have set $Z = V^*$. Hence in all cases (whether $P(G^*)$ is reached or not) BBC must output $Z_F = V^*$. □

We have shown that BBC converges to $V^*$. Previously, we alluded to the desirability of maximizing $f - \delta$ over $S$ by using algorithm BBC, rather than reformulation, if $f - \delta$ has certain attractive features. Presently we explicate the claim by demonstrating how LS may be adapted into a heuristic for solving $P(\hat{G})$ with a non-trivial, but not a priori bounded, approximation ratio if $f - \delta$ is non-negative.

A key observation is that both $2^{I(\hat{G})} \cap S$ and $\{\emptyset, v_1, ..., v_n\} \cap S$ are subset systems. Thus $\hat{F}$ is a subset system, and $P(\hat{G}) = \{\max (f_u - \hat{\delta})(V') \mid V' \in \hat{F}\}$ is equivalent to the unconstrained problem $\{\max (f_u - \hat{\delta} - Mq)(V') \mid V' \in 2^V\}$, where $M > \max_{V' \in 2^V} |f_u - \hat{\delta}|$ and $q$ is as described in the proof of Prop. 14, by Cor. 15. We assume the existence of an oracle which can determine membership in $S$ in polynomial time, and thus there exists a fixed polynomial $Q(n)$ bounding the complexity of evaluating $q$. Since $f_u - \hat{\delta}$ is submodular, so is $f_u - \hat{\delta} - Mq$, which is our starting point for adapting LS. We shall also assume that $(f - \delta)(v) > 0 \; \forall v \in V$, and hence $(f_u - \hat{\delta})(v) > 0 \; \forall v \in V$. This is done without loss of (too much) generality, as $f - \delta$ is non-negative and thus an arbitrarily small perturbation of $f - \delta$ on these n singletons is sufficient.

**Alg. LSA (Adapted)**

Input: A normalized, non-negative, submodular function g: $2^X \to \mathbb{R}$ with $g(x) > 0$ $\forall x \in X$, where X is a set of size r. A non-empty subset system $S$, whose membership function can be evaluated in polynomial time, and such that $\forall x \in X, x \in S$. $\varepsilon > 0$, and a constant M and a function q: $2^X \to \mathbb{R}$ defined as in Cor. 15.

Output: An approximate maximizer of g restricted to $S$.

0) Initialize $i = 0$, $X^{(0)} = X$, $X' = \emptyset$.

1) Let $X' = \text{argmax}_{x \in X^{(i)}} (g-Mq)(x)$.

2) If $\exists x \in X^{(i)} \setminus X'$ such that $(g-Mq)(X' \cup x) > (1 + \frac{\varepsilon}{r^2})(g-Mq)(X')$, update $X' \to X' \cup x$ and go back to step 2.

3) If $\exists x \in X'$ such that $(g-Mq)(X' \setminus x) > (1 + \frac{\varepsilon}{r^2})(g-Mq)(X')$, update $X' \to X' \setminus x$ and go back to step 2.

4) Output $\hat{X}^{(i+1)} = X'$, $\hat{X}^{(i+1)C} = X^{(i)} \setminus \hat{X}^{(i+1)}$.

    a) If $\hat{X}^{(i+1)C} \notin S$, set $i \to i+1$, $X^{(i+1)} \to \hat{X}^{(i+1)C}$, $r = |X^{(i+1)}|$, $X' \to \emptyset$, and go to step 1.

    b) Else, update $X^{(i+1)} \to \hat{X}^{(i+1)C}$, and call LS as a subroutine to maximize $(g-Mq)$ over $2^{X^{(i+1)}}$, and return approximate solutions $\hat{X}^{(i+2)}$ and $\hat{X}^{(i+2)C}$. Output $\hat{X}_A = \text{argmax}\{(g-Mq)(X') \mid X' = \hat{X}^{(i)} \text{ or } \hat{X}^{(i)C}, i = 0, 1,...\}$.

    Note that steps 1 through 3 of LSA constitute calling LS as a subroutine. In step 4, we differentiate between returning to step 1 and calling LS as a subroutine in order to make the convergence criteria for the algorithm clearer. We say that LSA has completed an iteration whenever step 4 is reached. Each iteration of LSA is comprised of several subiterations which are repetitions of steps 2 and 3, i.e. iterations of the LS subroutine.

    Before deriving the approximation ratio of LSA, we must prove that is converges. We say that LSA has completed an iteration whenever step 4 is reached. Each iteration of LSA is constituted by several subiterations of steps 2 and 3. Our approach to proving LSA converges will be to show that each iteration of LSA requires as most $O(\frac{1}{\varepsilon}Q(r)*r^3 \log r)$ operations, where Q(r) is a polynomial bounding the complexity of evaluation g-Mq. We then show that LSA requires finitely many, specifically, no greater than r, iterations before the stopping condition in step 4 is satisfied. Throughout the next several proofs, it will be helpful to recall that $(g-Mq)(X') \geq 0$ iff $X' \in S$.

**Prop. 17:** Each iteration of LSA requires at most $O(\frac{1}{\varepsilon}Q(r)*r^3 \log r)$ operations.

Proof: Our proof follows that of Thm. 3.4 of [4]. Consider the ith iteration of LSA. If $X^{(i)} = \varnothing$, the claim is trivially true. If not, then let $\hat{x}$ be an element of $X^{(i)}$ that maximizes (g-Mq) over singletons. Notice that $(g\text{-}Mq)(\hat{x}) = g(\hat{x}) > 0$. Since X' is initialized to $\varnothing$ and g is normalized, it follows that $(g\text{-}Mq)(X') = (g\text{-}Mq)(\varnothing) = g(\varnothing) = 0$. Thus $X' \to \hat{x}$ in step 1. Because g is submodular we know that $\max_{X' \in 2^{X^{(i)}}} g(X') \le rg(\hat{x})$, and so $\max_{X' \in 2^{X^{(i)}}} (g\text{-}Mq)(X') \le r(g\text{-}Mq)(\hat{x})$, where r denotes $|X^{(i)}|$. Since (g-Mq)(X') increases by a factor of $(1+\frac{\varepsilon}{r^2})$ in each subiteration of LSA, it follows that after k+1 subiterations, $(g\text{-}Mq)(X') \ge (1+\frac{\varepsilon}{r^2})^k (g\text{-}Mq)(\hat{x})$. Hence after at most $O(\frac{1}{\varepsilon} r^2 \log r)$ subiterations we must have $(g\text{-}Mq)(X') \ge r(g\text{-}Mq)(\hat{x})$ and so LSA finished its iteration. Since each step requires no more than Q(r)*r operations, the desired result is consequent.

**Prop. 18:** LSA requires at most r iterations.

Proof: Since $\varnothing \in S$, the algorithm clearly terminates if $\hat{X}^{(i+1)C} = \varnothing$ for some i. Therefore to prove the claim it is sufficient to show that $|\hat{X}^{(i+1)C}|$ is strictly decreasing in I, which is equivalent to showing that $\hat{X}^{(i+1)}$ is non-empty if $X^{(i)} \ne \varnothing$. The last claim follows by the proof of Prop. 17 which shows that $(g\text{-}Mq)(\hat{X}^{(i+1)}) > 0$, and thus $\hat{X}^{(i+1)} \ne \varnothing$, whenever $X^{(i)} \ne \varnothing$.

Hence LSA requires $O(\frac{1}{\varepsilon} Q(r) * r^4 \log r)$ operations. For convenience we let $r^{(i)}$ denote $|X^{(i)}|$ as we derive the approximation guarantee given by LSA. By design of LSA, during each iteration i, the resultant $\hat{X}^{(i+1)}$ is "locally optimal" in the sense that

$(1+\frac{\varepsilon}{(r^{(i)})^2})(g\text{-}Mq)(\hat{X}^{(i+1)}) > (g\text{-}Mq)(\hat{X}^{(i+1)} \cup x) \ \forall x \in X^{(i)} \setminus \hat{X}^{(i+1)}$, and

$(1+\frac{\varepsilon}{(r^{(i)})^2})(g\text{-}Mq)(\hat{X}^{(i+1)}) > (g\text{-}Mq)(\hat{X}^{(i+1)} \setminus x) \ \forall x \in \hat{X}^{(i+1)}$. N.b. that this relies implicitly on the fact that (g-Mq)(X') becomes, and stays, positive in step 1.

**Prop. 19**: Let $\hat{X}^{(i+1)}$ be locally optimal. Then for any $X' \in 2^{X^{(i)}}$ that is comparable to $\hat{X}^{(i+1)}$, $(1+\frac{\varepsilon}{r^{(i)}})(g\text{-}Mq)(\hat{X}^{(i+1)}) > (g\text{-}Mq)(X')$.

Proof: See [4] Lemma 3.3.

**Prop. 20**: Suppose that LSA converges at the end of iteration $i = k$, $k \in \{0,...,r-1\}$, then the solution $\hat{X}_A$ output by LSA satisfies

$[(k+4)+\varepsilon H(k+2)+8\varepsilon](g-Mq)(\hat{X}_A) > \max_{X' \subset X}(g-Mq)(X')$.

Proof: Let $\hat{X}_{i+1}$ denote an optimal solution to $\{\max (g-Mq)(X') \mid X' \subset X^{(i)}\}$, where $X^{(0)} = X$ and $X^{(i+1)} = \hat{X}^{(i)C}$. Thus $\hat{X}_1$ is an optimal solution to our original problem. Let $\hat{X}^{(i+1)}$ and $\hat{X}^{(i+1)C}$ denote the solution output by LSA, and its complement, for iteration $i = 0,1,...,k$. By the submodularity of $(g-Mq)$ we have
$(g-Mq)(\hat{X}_1 \cap \hat{X}^{(1)}) + (g-Mq)(\hat{X}_1 \setminus \hat{X}^{(1)}) \geq (g-Mq)(\hat{X}_1) + (g-Mq)(\emptyset)$. Since g and q are normalized, the latter because $\emptyset \in S$, this implies that
$(g-Mq)(\hat{X}_1 \cap \hat{X}^{(1)}) + (g-Mq)(\hat{X}_1 \setminus \hat{X}^{(1)}) \geq (g-Mq)(\hat{X}_1)$. Clearly $\hat{X}_1 \cap \hat{X}^{(1)}$ is comparable to $\hat{X}^{(1)}$ and so $(1 + \frac{\varepsilon}{r^{(0)}})(g-Mq)(\hat{X}^{(1)}) > (g-Mq)(\hat{X}_1 \cap \hat{X}^{(1)})$.

Furthermore, $(g-Mq)(\hat{X}_1 \setminus \hat{X}^{(1)}) \leq \max_{X' \subset X^{(1)}}(g-Mq)(X')$ as $\hat{X}_1 \setminus \hat{X}^{(1)} \subset \hat{X}^{(1)C} = X^{(1)}$.
Iteration $i = 1$ of LSA attempts to maximize $(g-Mq)$ over $2^{X^{(1)}}$ and outputs an approximate maximizer $\hat{X}^{(2)}$. Applying submodularity again, we have
$(1 + \frac{\varepsilon}{r^{(1)}})(g-Mq)(\hat{X}^{(2)}) + (g-Mq)(\hat{X}_2 \setminus \hat{X}^{(2)}) > (g-Mq)(\hat{X}_2)$, and thus
$(1 + \frac{\varepsilon}{r^{(0)}})(g-Mq)(\hat{X}^{(1)}) + (1 + \frac{\varepsilon}{r^{(1)}})(g-Mq)(\hat{X}^{(2)}) + (g-Mq)(\hat{X}_2 \setminus \hat{X}^{(2)}) > (g-Mq)(\hat{X}_1)$. Repeating this procedure k times gives us
$(1 + \frac{\varepsilon}{r^{(0)}})(g-Mq)(\hat{X}^{(1)}) + ... + (1 + \frac{\varepsilon}{r^{(k)}})(g-Mq)(\hat{X}^{(k+1)}) + (g-Mq)(\hat{X}_{k+1} \setminus \hat{X}^{(k+1)}) > (g-Mq)(\hat{X}_1)$.

By assumption, the algorithm converges at the end of iteration k, therefore we may assume that $|X^{(k)}| \geq 2$, otherwise we would have converged in an earlier iteration, as all singletons and $\emptyset$ are members of $S$, and that $X^{(k+1)} \in S$. Since $X^{(k+1)} \in S$, $(g-Mq)$ is non-negative on $2^{X^{(k+1)}}$, and so the LS subroutine of LSA returns $\hat{X}^{(k+2)}$ and $\hat{X}^{(k+2)C}$, one of which is a $(\frac{1}{3} - \frac{\varepsilon}{r^{(k+1)}})$-approximate maximizer of $(g-Mq)$ over $2^{X^{(k+1)}}$. Without loss of generality, $(g-Mq)(\hat{X}^{(k+2)}) \geq (g-Mq)(\hat{X}^{(k+2)C})$, which gives us
$(\frac{3r^{(k+1)}}{r^{(k+1)}-3\varepsilon})(g-Mq)(\hat{X}^{(k+2)}) \geq \max_{X' \subset X^{(k+1)}}(g-Mq)(X')$, and so
$(3+9\varepsilon)(g-Mq)(\hat{X}^{(k+2)}) \geq \max_{X' \subset X^{(k+1)}}(g-Mq)(X')$. In turn, we have
$(1+\frac{\varepsilon}{r^{(0)}})(g-Mq)(\hat{X}^{(1)}) + ... + (1+\frac{\varepsilon}{r^{(k)}})(g-Mq)(\hat{X}^{(k+1)}) + (3+9\varepsilon)(g-Mq)(\hat{X}^{(k+2)}) > (g-Mq)(\hat{X}_1)$.

Therefore $[(k+4) + \frac{\varepsilon}{r^{(0)}} + ... + \frac{\varepsilon}{r^{(k)}} + 9\varepsilon](g-Mq)(\hat{X}_A) > (g-Mq)(\hat{X}_1)$.

Finally, since $r^{(0)} > r^{(1)} > \ldots > r^{(k)} \geq 2$, by the convergence proof of Prop. 18 and our previous observation on $|X^{(k)}|$, we see that $\sum_{i=0}^{k} \frac{\varepsilon}{r^{(i)}} \leq \varepsilon(H(k+2)-1)$, where $H(k)$ is the kth harmonic number. Hence $[(k+4)+\varepsilon H(k+2)+8\varepsilon](g\text{-}Mq)(\hat{X}_A) > (g\text{-}Mq)(\hat{X}_1)$.

Substituting $I(\hat{G})$ for X, and $f_u - \hat{\delta}$ for g, yields a heuristic procedure for solving $P(\hat{G})$ for $\hat{G} \neq K_n$. (For $\hat{G}=K_n$ the problem is trivial to solve, only requiring evaluating the function on at most every singleton and $\varnothing$.) A suitable choice of upper bound M is given by $|I(\hat{G})|*(f_u-\hat{\delta})(\hat{v})$, where $\hat{v}$ maximizes $f_u - \hat{\delta}$ over all singletons in $I(\hat{G})$. The assumption $v \in S \; \forall v \in I(\hat{G})$ is valid, for if not we could have reduced the search space of solutions to $2^{V \cap S}$ at the outset of algorithm BBC. Finally, we observe that BBC can be adapted to make use of approximate solutions to $P(\hat{G})$ in a fashion similar to the adaptation of BB.

## §6. Conclusions

Our primary contribution has been to demonstrate a "universal decomposition" of any set function $\theta: 2^V \to \mathbb{R}$ as $f - \delta$, and to use this decomposition to maximize $\theta$ via a branch and bound framework. Since most of our results follow from this, we briefly compare our method, BB, to the best-known method for maximizing the difference of two submodular functions, the Submodular-Supermodular Procedure of [10]. Note that SSP permits the maximization of the difference of two arbitrary submodular functions f and f', whereas BB requires that $f' = \delta$ for some graph G. By Prop. 3 this difference is largely cosmetic, as $f - f'$ can be reformulated as the difference of a submodular function and a cut function, though finding a good reformulation is non-trivial.

The two algorithms employ fundamentally different approaches. SSP is an ascent method, producing a monotonically improving sequence of intermediate solutions, and converges to a locally-optimal solution. BB, on the other hand, guarantees global optimality but does not necessarily find a better solution after solving each subproblem. Both methods have exponential worst case complexity. The claim has already been shown for BB, but is not immediately obvious for SSP as the latter guarantees only *local* optimality. However, during each iteration of SSP, with current best solution V' of size k, one must select a permutation $\pi$ of the ground set V so that $\{v_{\pi(1)},\ldots,v_{\pi(k)}\}=V'$ and solve a (polynomial-time) maximization problem dependent upon $\pi$. It can be shown [9] that if $\pi$ is selected arbitrarily from amongst all permutations satisfying $\{v_{\pi(1)},\ldots,v_{\pi(k)}\}=V'$, that SSP requires $O(2^n)$ iterations to converge. Additionally, SSP does not appear to be interruptible.

The ascent property of SSP suggests that it may return superior solutions during the early stages of the algorithm, but with no guaranteed (additive or multiplicative) difference from optimality. Because of the vastly different approaches taken by the algorithms, extensive numerical testing is warranted to benchmark the relative performance of BB versus SSP, and to identify classes of problems amenable to either or both approaches. These results will be communicated in a subsequent paper.